\documentclass{amsart}
\usepackage[all]{xy}

\DeclareMathSymbol{\twoheadrightarrow}  {\mathrel}{AMSa}{"10}

\def\Q{{\mathbb Q}}
\def\Z{{\mathbb Z}}
\def\C{{\mathbb C}}

\def\P{{\mathbb P}}

\def\Gal{\mathrm{Gal}}

\def\End{\mathrm{End}}
\def\Aut{\mathrm{Aut}}
\def\Hom{\mathrm{Hom}}

\def\L{\mathcal{L}}

\def\I{{\mathcal I}}

\def\fchar{\mathrm{char}}

    \def\G{\mathbf{G}}

       \def\cl{\mathrm{cl}}
       \def\Spec{\mathrm{Spec}}

\def\pr{\mathrm{pr}}
\def\M{\mathcal{M}}
    \def\m{\mathrm{m}}

\def\dim{\mathrm{dim}}
                     \def\Mat{\mathrm{Mat}}

\def\OC{{\mathcal O}}

\def\q{{\mathfrak q}}

\def\b{{\mathfrak b}}

\def\m{{\mathfrak m}}

          \def\rk{\mathrm{rk}}
          \def\j{\mathrm{j}}
\newtheorem{thm}{Theorem}[section]
\newtheorem{lem}[thm]{Lemma}
\newtheorem{cor}[thm]{Corollary}
\newtheorem{prop}[thm]{Proposition}
\theoremstyle{definition}
\newtheorem{defn}[thm]{Definition}
\newtheorem{ex}[thm]{Example}
\newtheorem{rem}[thm]{Remark}

            \newtheorem{sect}[thm]{}
\title[Homomorphisms of abelian varieties over finite fields]
{Homomorphisms of abelian varieties over finite fields}
\author[Yuri\ G.\ Zarhin]{Yuri\ G.\ Zarhin}
\address{Department of Mathematics, Pennsylvania State University,
University Park, PA 16802, USA}
\address{Institute for Mathematical Problems in Biology of the Russian Academy of
Sciences, Pushchino, Moscow Region, Russia}
 \email{zarhin\char`\@math.psu.edu}

\begin{document}




\maketitle

The aim of this note is to give a proof of Tate's theorems on homomorphisms of
abelian varieties over finite fields and the corresponding $\ell$-divisible
groups \cite{Tate, MW}, using ideas of \cite{ZarhinMatZam,ZarhinJ}. We give a
unified treatment for both $\ell\ne p$ and $\ell=p$ cases. In fact, we prove a
slightly stronger version of those theorems with ``finite coefficients". We use
neither  the existence (and properties) of the Frobenius endomorphism (for
$\ell\ne p$) nor Dieudonne modules (for $\ell=p$).

The paper is organized as follows. (A rather long) Section \ref{one} contains
auxiliary results about finite commutative group schemes and abelian varieties
with special reference to isogenies and polarizations. We discuss
$\ell$-divisible groups (aka Barsotti--Tate groups) in Section \ref{ellDIV}.
Section \ref{use} contains useful results that play a crucial role in the proof
of main results that are stated in Section \ref{mainres}.

The next five Sections  contain proofs of  results that were stated in Section
\ref{use}. In Section \ref{pLOZ} we discuss abelian subvarieties of a given
abelian variety. Section \ref{polAV} deals with the finiteness of the set of
abelian varieties of given dimension and ``bounded degree" over a finite field.
In Section \ref{quatTRICK} we present a so called {\sl quaternion trick}. In
Section \ref{finiteKER} we prove a crucial result about arbitrary finite group
subschemes of abelian varieties over finite fields. In Section \ref{finiteDIV}
we try to divide endomorphisms of a given abelian variety modulo $n$.

 The main results of this paper are proven in
  Section \ref{proofMAIN}. Their variants for Tate modules are
 discussed in Section \ref{TateMOD}. An example of non-isomorphic elliptic
 curves over a finite field with isomorphic $\ell$-divisible groups (for all
 primes $\ell$) is discussed in Section \ref{nonISOM}.

I am grateful to Frans Oort and Bill Waterhouse for useful discussions and to
the referee, whose comments helped to improve the exposition. My special thanks
go to Dr. Boris Veytsman for his help with \TeX nical problems.

\section{Definitions and statements}
\label{one}

Throughout this paper $K$ is a field and $\bar{K}$ its algebraic
closure. If $X$ (resp. $W$) is an algebraic variety (resp. group
scheme) over $K$ then we write $\bar{X}$ (resp. $\bar{W}$) for the
corresponding algebraic variety $X\times_{\Spec(K)}\Spec(\bar{K})$
(resp. group scheme $W\times_{\Spec(K)}\Spec(\bar{K})$) over
$\bar{K}$. If $f:X\to Y$ is a a regular map of algebraic varieties
over $K$ then we write $\bar{f}$ for the corresponding map
$\bar{X}\to \bar{Y}$.

\begin{sect} {\bf Finite commutative group schemes over fields}.
We refer the reader to the books of Oort \cite{Oort}, Waterhouse
\cite{W} and Demazure--Gabriel \cite{DG} for basic properties of
commutative group schemes; see also \cite{Shatz,Pink}.

Recall that a group scheme $V$ over $K$ is called finite if the
structure morphism $V \to \Spec(K)$ is finite. Since $\Spec(K)$ is
a one-point set, it follows from the definition of finite morphism
\cite[Ch. II, Sect. 3]{Ha} that $V$ is an affine scheme and
$\Gamma(V,\OC_V)$ is a finite-dimensional commutative $K$-algebra.
 The $K$-dimension of the $\Gamma(V,\OC_V)$ is called
 the {\sl order}  of $V$ and  denoted by $\#(V)$. An analogue of
 Lagrange theorem \cite{OT} asserts that multiplication by $\#(V)$
 kills commutative $V$.

 Let $V$ and $W$ be finite commutative group schemes over $K$ and let
$u:V \to W$ be a morphism of group $K$-schemes. Both  $V$ and $W$
are affine schemes, $A=\Gamma(V,\OC_V)$ and $B=\Gamma(W,\OC_W)$
are finite-dimensional (commutative) $K$-algebras (with $1$),
$V=\Spec(A), W=\Spec(B)$ and $u$ is induced by a certain
$K$-algebra homomorphism
$$u^{*}:B \to A.$$
Since $V$ and $W$ are commutative group schemes, $A$ and $B$ are
cocommutative  Hopf $K$-algebras. Since $u$ is a morphism of group
schemes, $u^{*}$ is a morphism of Hopf algebras. It follows that
$C:=u^{*}(B)$ is a $K$-subalgebra and also a Hopf subalgebra in
$A$. It follows that $U:=\Spec(C)$ carries the natural structure
of a finite group scheme over $K$ such that the natural scheme
morphism $U\to V$ induced by $u^{*}: B \twoheadrightarrow
u^{*}(B)=C$ is a morphism of group schemes. In addition, the
inclusion $C\subset A$ induces the morphism of schemes $V \to U$,
which is also a morphism of group schemes. The latter morphism is
an epimorphism in the category of finite commutative group schemes
over $K$, because the corresponding map
$$C=\Gamma(U,\OC_U)\to \Gamma(V,\OC_V)=A$$
is nothing else but the inclusion map $C\subset A$ and therefore
is injective \cite{OS} (see also \cite{HM}).

On the other hand, the surjection $B\twoheadrightarrow C$ provides
us with a canonical isomorphism $U \cong \Spec(B/\ker(u^{*}))$; in
addition, we observe that $\Spec(B/\ker(u^{*}))$ is a (closed)
group subscheme of $\Spec(B)=W$. We denote $\Spec(B/\ker(u^{*}))$
by $u(V)$ and call it the image of $u$ or the image of $V$ with
respect to $u$ and denote by $u(V)$. Notice that the set theoretic
image of $u$ is closed and our definition of the image of $u$
coincides with the one given in \cite[Sect. 5.1.1]{EH}.

One may easily check that the closed embedding
$j:u(V)\hookrightarrow V$ induced by $B\twoheadrightarrow
B/\ker(u^{*})$ is an image in the category of (affine) schemes
over $K$. This means that if $\alpha, \beta: W \to S$  are two
morphisms of schemes over $K$ such that their {\sl restrictions}
to $u(V)$ do coincide, i.e., $\alpha j=\beta j$ (as morphisms from
$u(V)$ to $S$) then $\alpha u=\beta u$ (as morphisms from $U$ to
$S$). It follows that $j$ is also an image in the category of
finite commutative group schemes. group \cite[Sect. 10]{Pink}.
\end{sect}

\begin{thm}[Theorem of Gabriel \cite{OS,HM}]
\label{gabriel}
 The category of finite commutative group schemes
over a field is abelian.
\end{thm}

\begin{rem}
\label{orderE} Let $V$ be a finite commutative group scheme over
$K$ and let $W$ be its finite closed group subscheme. If $V \to U$
is a {\sl surjective} morphism of finite commutative group schemes
over $K$ with kernel $W$ then \cite{HM}
$$\#(V)=\#(W)\cdot \#(U).$$
Recall that $\Gamma(W,\OC_W)$ is the quotient of
$\Gamma(V,\OC_V)$. In particular,
 if the orders of $V$ and $W$ do
coincide then $V=W$.
\end{rem}

\begin{sect} {\bf Abelian varieties over fields}.
 We refer the reader to the books of Mumford \cite{MumfordAV},
 Shimura \cite{Shimura}   for basic properties of abelian varieties
 (see also Lang's book  \cite{Lang} and papers of Waterhouse \cite{W2}, Deligne \cite{Deligne}, Milne \cite{Milne} and Oort \cite{OG}).
If  $X$ is an abelian variety  over $K$ then we write $\End(X)$ for the ring of
all $K$-endomorphisms of $X$. If $m$ is an integer then write $m_X$ for the
multiplication by $m$ in $X$; in particular, $1_X$ is the identity map.
(Sometimes we will use notation $m$ instead of $m_X$.)

If $Y$ is an abelian variety  over $K$ then we write $\Hom(X,Y)$ for the group
of all $K$-endomorphisms $X\to Y$.

\begin{rem}
Warning: sometimes in the literature, including my own papers, the notation
$\End(X)$ is  used for the ring of $\bar{K}$-endomorphisms.
\end{rem}

It is well known \cite[Sect. 19, Theorem 3]{MumfordAV} that
$\Hom(X,Y)$ is a free commutative group of finite rank. We write
$X^t$ for the dual of $X$ (See \cite[ Sect. 9--10]{Milne} for the
definition and basic properties of the dual of an abelian variety.)
In particular, $X^t$ is also an abelian variety over $K$ that is
isogenous to $X$ (over $K$). If $u\in \Hom(X,Y)$ then we write $u^t$
for its dual in $\Hom(Y^t,X^t)$. We have
$$\bar{X}^t=\overline{X^t}.$$

If $n$ is a positive integer then we write $X_n$ for the kernel of
 $n_X$; it is a finite commutative
(sub)group scheme  (of $X$) over $K$ of rank $2\dim(X)$. By
definition, $X_n(\bar{K})$ is the kernel of multiplication by $n$
in $X(\bar{K})$.

If $n$ is not divisible by $\fchar(K)$ then $X_n$ is an \'etale
group scheme and  it is well-known \cite[Sect. 4]{MumfordAV} that
$X_n(\bar{K})$ is a free $\Z/n\Z$-module of rank $2\dim(X)$ and
all $\bar{K}$-points of $X_n$ are defined over a finite separable
extension of $K$. In particular, $X_n(\bar{K})$ carries a natural
structure of Galois module.
\end{sect}

\begin{sect} {\bf Isogenies}.
Let $W\subset X$ be a finite group subscheme over $K$. It follows
from the analogue of Lagrange theorem that $W\subset X_d$ for
$d=\#(W)$. The quotient $Y:=X/W$ is an abelian variety over $K$
and the canonical isogeny $\pi: X\to X/W=Y$ has kernel $W$ and
degree $\#(W)$ (\cite[Sect. 12, Corollary 1 to Theorem
1]{MumfordAV}, \cite[Sect.  2, pp. 307-314]{DG}). In particular,
every homomorphism of abelian varieties $u:X \to Z$ over $K$ with
$W\subset \ker(u)$ factors through $\pi$, i.e., there exists a
unique homomorphism of abelian varieties $v:Y\to Z$ over $K$ such
that
$$u=v\pi .$$
If $m$ is a positive integer then
$$\pi m_X=m_Y \pi \in \Hom(X,Y).$$
Let us put
 $$m^{-1}W:=\ker(\pi m_X)=\ker(m_Y \pi) \subset X.$$
 For every commutative $K$-algebra $R$ the group of $R$-points
 $m^{-1}W(R)$ is the set of all $x \in X(R)$ with
 $$m x\in W(R)\subset X(R).$$
 For example, if $W=X_n$ then $$Y=X, \pi=n_X, m^{-1}X_n=X_{nm}.$$
 In general, if $W\subset X_n$ then $m^{-1}W$ is a closed group subscheme in
 $X_{nm}$. E.g., $W$ is always a closed group subscheme of $X_{dm}$
 and therefore is a finite group subscheme of $X$ over $K$. The
 order
 $$\#(m^{-1}W)=\deg(\pi m_X)=\deg(\pi)\deg(m_X)=\#(W)\cdot m^{2\dim(X)}.$$
 We have
 $$X_m \subset m^{-1}W, \ m_X (m^{-1}W)\subset W$$
 and the kernel of $m_X:m^{-1}W \to W$ coincides with $X_m$.
\begin{lem}
\label{imageAG}
 The image
$ m_X (m^{-1}W)=W$.
\end{lem}
\begin{proof}
Let us denote the image by $G$. By Remark \ref{orderE}, $\#(G)$ is
the ratio $$\#(m^{-1}W)/\#(X_m)=\dim(W),$$ i.e., the orders of $G$
and $W$ do coincide. Since $G\subset W$, we have (by the same
Remark) $G=W$.
\end{proof}

\begin{ex}
\label{Xn} If $W=X_n$ then $m^{-1}X_n=X_{nm}$ and therefore
$m(X_{nm})=X_n$.
\end{ex}

\begin{lem}
\label{Xn1} If $r$ is a positive integer then  $r(X_n)=X_{n_1}$
where $n_1=n/(n,r)$.
\end{lem}

\begin{proof}
We have $r=(n,r)\cdot r_1$ where $r_1$ is a positive integer such
that $n_1$ and $r_1$ are relatively prime. This implies that
$r_1(X_{n_1})=X_{n_1}$. By Lemma \ref{Xn1}, $(n,r)(X_n)=X_{n_1}$.
This implies that
$$r(X_n)= r_1(n,r)(X_n)=r_1((n,r)(X_n))=r_1(X_{n_1})=X_{n_1}.$$
\end{proof}

\begin{lem}
\label{isogprime} Let $X$ and $Y$ be abelian varieties over a field $K$. Let
$u:X\to Y$ be a $K$-homomorphism of abelian varieties. Let $n>1$ be an integer
and $u_n:X_n \to Y_n$ the morphism of commutative group schemes over $K$
induced by $u$.

\begin{itemize}
\item[(i)]
 Suppose that $u$ is an isogeny and $\deg(u)$ and $n$ are
relatively prime. Then  $u_n:X_n \to Y_n$ is an isomorphism.

\item[(ii)] Suppose that $u_n:X_n \to Y_n$ is an isomorphism. Then $u$ is an
isogeny and $\deg(u)$ and $n$ are relatively prime.
\end{itemize}
\end{lem}

\begin{proof}
Let $u$ be an isogeny such that $m:=\deg(u)$ and $n$ are relatively prime. Then
 $\ker(u)\subset X_m$. It follows that there exists a $K$-isogeny $v:Y\to X$
 such that
 $$vu=m_X, uv=m_Y.$$

 {\bf (i)}. Since multiplication by $m$ is an automorphism of both $X_n$ and $Y_m$, we
 conclude that $u_n:X_n \to Y_n$ and $v_n:Y_n \to X_n$ are isomorphisms.

 {\bf (ii)}. Suppose that $u_n$ is an isomorphism. This implies that the
 orders of $X_n$ and $Y_n$ coincide and therefore $\dim(X)=\dim(Y)$. We need to
 prove that $u$ is isogeny and $\deg(u)$ and $n$ are relatively prime. In order
 to do that, we may assume that $K$ is algebraically closed (replacing $K,X,Y,u$ by
 $\bar{K},\bar{X},\bar{Y},\bar{u}$  respectively). Let us put $Z:=u(Y)\subset
 X$: clearly, $Z$ is a (closed) abelian subvariety of $Y$ and therefore $\dim(Z)\le \dim(Y)$.
 It is also clear that
 $u:X \to Y$ coincides with the composition of the natural surjection
  $X \to u(X)=Z$ and the inclusion map $j:Z \hookrightarrow X$. This implies
  that $u_n(X_n)$ is a (closed) group subscheme of $j_n(Z_n)\subset Y_n$. It
  follows that
  $$\#(u_n(X_n))\le \#(j_n(Z_n)) \le \#(Z_n)=n^{2\dim(Z)}.$$
  Since $u_n$ is an isomorphism, $u_n(X_n)=Y_n$ and therefore
$$\#(u_n(X_n))=\#(Y_n)=n^{2\dim(Y)}.$$
It follows that $$n^{2\dim(Y)}\le n^{2\dim(Z)}$$ and therefore $\dim(Y)\le
\dim(Z)$. (Here we use that $n>1$.) Since $Z$ is a closed subvariety in $Y$, we
conclude that $\dim(Z)=\dim(Y)$ and $Y=Z$. In other words, $u$ is surjective.
Taking into account that $\dim(X)=\dim(Y)$, we conclude that $u$ ia an isogeny.

Now let $m=d r$ where $d$ is the largest common divisor of $n$ and $m$. Then
$r$ and $n$ are relatively prime; in particular, multiplication by $r$ is an
automorphism of $X_n$. Let us denote $\ker(u)$ by $W$: it is a finite
commutative group scheme over $K$ of order $m$ and therefore
$$W \subset X_m.$$
This implies that for every commutative $K$-algebra $R$ we have
$$m \cdot W(R)=\{0\}.$$
On the other hand,  since $u_n$ is an isomorphism,  the kernel of $W(R)
\stackrel{n}{\to} W(R)$ is $\{0\}$. Since $d\mid n$, the kernel of $W(R)
\stackrel{d}{\to} W(R)$ is also $\{0\}$. This implies that $r\cdot W(R)=\{0\}$
for all $R$. Hence $W\subset X_r$. It follows that $\deg(u)=\#(W)$ divides
$\#(X_r)=r^{2\dim(X)}$ and therefore is coprime to $n$.
\end{proof}

The next statement will be used only in Section \ref{nonISOM}.

\begin{prop}
\label{composite} Let  $X$ and $Y$ be abelian varieties over a field $K$.
Suppose that for every prime $\ell$ there exists an isogeny $X\to Y$, whose
degree is not divisible by $\ell$. Then for every positive integer $n$ there
exists an isogeny $X\to Y$, whose degree is coprime to $n$. In particular, $X_n
\cong Y_n$.
\end{prop}

\begin{proof}
 Recall that the additive group $\Hom(X,Y)$ is isomorphic to
$\Z^{\rho}$ for some nonnegative integer $\rho$. In our case, $X$ and $Y$ are
isogenous over $K$ and therefore $\rho>0$.

Let $n$ be a positive integer and let $P(n)$ be the (finite) set of
its prime divisors. For each $\ell \in P(n)$ pick an isogeny
$v^{(\ell)}:X \to Y$, whose degree is not divisible by $\ell$. By
Lemma \ref{isogprime}(i), $v^{(\ell)}$ induces an isomorphism
$X_{\ell} \cong Y_{\ell}$.  Now, by the Chinese Remainder Theorem,
there exists $u \in \Hom(X,Y)\cong\Z^{\rho}$ such that
$$u -v^{(\ell)} \in \ell \cdot \Hom(X,Y) \ \forall \ \ell \in P.$$
This implies that for each $\ell \in P$ the homomorphisms $u$ and $v^{(\ell)}$
induce the same morphism $X_{\ell} \cong Y_{\ell}$, which, as we know, is an
isomorphism. It follows from  Lemma \ref{isogprime}(ii) that $u$ is an
isogeny, whose degree is not divisible by $\ell$. Hence $\deg(u)$ and $n$ are
coprime. Applying again Lemma \ref{isogprime}(i), we conclude that $u$ induces
an isomorphism $X_n \cong Y_n$.
\end{proof}

\end{sect}

\begin{sect} {\bf Polarizations}.
 A homomorphism $\lambda: X \to X^t$
is a {\sl polarization} if there exists an ample invertible sheaf
$\L$ on $\bar{X}$ such that $\bar{\lambda}$ coincides with
$$\Lambda_{\L}:\bar{X}^t\to \bar{X}^t, \ z \mapsto \cl(T_{z}^{*}L\otimes
L^{-1})$$ where $T_z: \bar{X}\to \bar{X}$ is the translation map
$$x \mapsto x+z$$
and $\cl$ stands for the isomorphism class of an invertible sheaf.
Recall \cite[Sect, 6, Proposition 1; Sect. 8, Theorem 1; Sect. 13,
Corollary 5]{MumfordAV} that a polarization is an {\sl isogeny}. If
$\lambda$ is an isomorphism, i.e., $\deg(\lambda)=1$, we call
$\lambda$ a {\sl principal polarization} and the pair $(X,\lambda)$
is called a principally polarized abelian variety (over $K$).

If $n:=\deg(\lambda)=\#(\ker(\lambda))$ then $\ker(\lambda)$ is
killed by multiplication by $n$, i.e., $\ker(\lambda)\subset X_n$.
For every positive integer $m$ we write $\lambda^m$ for the
polarization

$$X^m \to (X^m)^t= (X^t)^m, \ (x_1, \dots , x_m)\mapsto (\lambda(x_1), \dots ,
\lambda(x_m))$$ that corresponds to the ample invertible sheaf $\otimes_{i=1}^m
\pr_i^{*}\L$ where $\pr_i: X^m \to X$ is the $i$th projection map. We have
$$\dim(X^m)=m\cdot \dim(X), \ \deg(\lambda^m)=\deg(\lambda)^m$$
and $\ker(\lambda^m)={\ker(\lambda)}^m \subset (X^m)_n$ if
$\ker(\lambda)\subset X_n$.

There exists a {\sl Riemann form} - a skew-symmetric pairing of
group schemes over $\bar{K}$ \cite[Sect. 23]{MumfordAV}

$$e_{\lambda}:\ker(\bar{\lambda})\times \ker(\bar{\lambda}) \to \G_{\m}$$
where $\G_{\m}$ is the multiplicative group scheme over
$\bar{K}$.

If
$$e_{\lambda^m}:\ker(\bar{\lambda}^m) \times \ker(\bar{\lambda}^m) \to \G_{\m}$$
is the Riemann form for $\lambda^{m}$ then in obvious notation
$$e_{\lambda^m}(x,y)=\prod_{i=1}^m e_{\lambda}(x_i,y_i)$$
where $$x=(x_1,\dots , x_m), \ y=(y_1,\dots , y_m) \in
\ker(\bar{\lambda})^m = \ker(\bar{\lambda}^m).$$

We have

$$\Mat_m(\Z)\subset \Mat_m(\End(X))=\End(X^m).$$

One may easily check that every $u\in \Mat_m(\Z)$ leaves the group
subscheme $\ker(\bar{\lambda}^m)$ invariant and
$$e_{\lambda^m}(ux,y)=e_{\lambda^m}(x,u^{*} y)$$
where $u^{*}$ is the transpose of the matrix $u$. Notice that
$u^{*}$ viewed as an element of $$\Mat_m(\Z)\subset
\Mat_m(\End(X^t))=\End((X^t)^m)$$ coincides with  $u^t\in
\End((X^m)^{t})$.
\end{sect}

\begin{sect} {\bf Polarizations and isogenies}.
 Let $W\subset \ker(\lambda)$ be a finite group subscheme
over $K$. Recall that $Y:=X/W$ is an abelian variety over $K$ and
the canonical isogeny $\pi: X\to X/W=Y$ has kernel $W$ and degree
$\#(W)$.

 Suppose that $\bar{W}$ is isotropic with respect to
$e_{\lambda}$, i.e., the restriction of $e_{\lambda}$ to $\bar{W}\times
\bar{W}$ is trivial. Then there exists an ample invertible sheaf $\M$ on
$\bar{Y}$ such that $\L\cong \bar{\pi}^* \bar{\M}$ \cite[Sect. 23, Corollary to
Theorem 2, p. 231]{MumfordAV} and the $\bar{K}$-polarization
$\Lambda_{\bar{\M}}:\bar{Y}\to \bar{Y^t}$ satisfies
$$\bar{\lambda}= \overline{\pi^t} \Lambda_{\bar{\M}}\bar{\pi}.$$
Since $\bar{\pi^t}$ and $\bar{\pi}$ are isogenies that are defined
over $K$, the polarization $\Lambda_{\bar{\M}}$ is also defined
over $K$, i.e., there exists a $K$-isogeny $\mu:Y\to Y^t$ such
that $\Lambda_{\bar{\M}}=\bar{\mu}$ and
$$\lambda=\pi^t \mu \pi .$$
It follows that
$$\deg(\lambda)=\deg(\pi)\deg(\mu)\deg(\pi^t)=\deg(\pi)^2\deg(\mu)=(\#(W))^2\deg(\mu).$$
Therefore $\mu$ is a principal polarization (i.e., $\deg(\mu)=1$)
if and only if
$$\deg(\lambda)=(\#(W))^2.$$
\end{sect}

\section{$\ell$-divisible groups,  abelian varieties  and Tate modules}
 \label{ellDIV}
 Let $h$ be a non-negative integer and $\ell$  a prime. The
following notion was introduced by Tate \cite{TateP,Shatz}.

\begin{defn}
A $\ell$-divisible group $G$ over $K$ of height $h$ is a sequence
$\{G_{\nu},i_\nu\}_{\nu=1}^{\infty}$ in which:
\begin{itemize}
\item $G_{\nu}$ is a finite commutative group scheme over $K$ of order
$\ell^{h\nu}$.
 \item
 $i_{\nu}$ is a closed embedding $G_\nu\hookrightarrow G_{\nu+1}$ that is
 a morphism of group schemes. In addition, $i_{\nu}(G_{\nu})$ is the
 kernel of multiplication by $\ell^{\nu}$ in $G_{\nu+1}$.
 \end{itemize}
\end{defn}

\begin{ex}
Let $X$ be an abelian variety over $K$ of dimension $d$. Then it is known
\cite{TateP,Shatz} that the sequence $\{X_{\ell^{\nu}}\} _{\nu=1}^{\infty}$ is
an $\ell$-divisible group over $K$ of height $2d$. Here $i_{\nu}$ is the {\sl
inclusion map} $X_{\ell^{\nu}}\hookrightarrow X_{\ell^{\nu+1}}$. We denote this
$\ell$-divisible group by $X(\ell)$.
\end{ex}

\begin{sect}{\bf Homomorphisms of $\ell$-divisible groups and abelian
varieties}. If $H=\{H_{\nu},j_{\nu}\}_{\nu=1}^{\infty}$ is an $\ell$-divisible
group  over $K$ then a morphism $u: G \to H$ is a sequence
$\{u_{(\nu)}\}_{\nu=1}^{\infty}$ of morphisms of group schemes over $K$
$$u_{(\nu)}: G_{\nu} \to H_{\nu}$$
such that the composition
$$u_{(\nu+1)}i_{\nu}: G_{\nu} \hookrightarrow G_{\nu+1} \to H_{\nu+1}$$
coincides with
$$j_{\nu}u_{(\nu)}:G_{\nu} \to H_{\nu} \hookrightarrow H_{\nu+1},$$
i.e., the diagram

\begin{displaymath}
 \xymatrix{
   G_{\nu} \ar[r]^{u_{(\nu)}} \ar[d]_{i_{\nu}}  & H_{\nu} \ar[d]^{j_{\nu}} \\
   G_{\nu+1} \ar[r]^{u_{(\nu+1)}}                   & H_{\nu+1}
}
\end{displaymath}
is commutative.

\begin{rem}
A morphism $u$ is an isomorphism of $\ell$-divisible groups if and only if all
$u_{(\nu)}$ are isomorphisms of the corresponding finite group schemes.
\end{rem}

 The group $\Hom(G,H)$ of morphisms from $G$ to $H$ carries a
natural structure of $\Z_{\ell}$-module induced by the natural structures of
$\Z/\ell^v=\Z_{\ell}/\ell^{\nu}$-module on $\Hom(G_{\nu},H_{\nu})$. Namely, if
$u=\{u_{(\nu})\}_{\nu=1}^{\infty}\in \Hom(G,H)$ and $a \in \Z_{\ell}$ then
$au=\{(au)_{(\nu)}\}_{\nu=1}^{\infty}$ may be defined as follows. For each
$\nu$ pick $a_{\nu} \in \Z$ with $a-a_{\nu}\in \ell^{\nu} \Z_{\ell}$ and put
$$(au)_{(\nu)}:=a_{\nu} u_{(\nu)}: G_{\nu} \to H_{\nu}.$$
Since multiplication by $\ell^{\nu}$ kills $G_{\nu}$, the definition of
$(au)_{(\nu)}$ does not depend on the choice of $a_{\nu}$.

 Let $X$ and $Y$ be abelian varieties over $K$. There is a natural homomorphism of commutative groups
 $\Hom(X,Y) \to \Hom(X(\ell),Y(\ell))$. Namely, if $u \in
 \Hom(X,Y)$ then $u(X_{\ell^{\nu}})$ lies in the kernel of
 multiplication by $\ell^{\nu}$, i.e.
 $u(X_{\ell^{\nu}})\subset Y_{\ell^{\nu}}$. In fact, we get the natural
 homomorphism
 $$\Hom(X,Y)\otimes \Z/\ell^{\nu} \to \Hom(X_{\ell^{\nu}},Y_{\ell^{\nu}}),$$
 which is known to be an embedding. (See also Lemma \ref{divprin}
 below.)

 Since
 $\Hom(X(\ell),Y(\ell))$ is a $Z_{\ell}$-module, we get the
 natural homomorphism of $\Z_{\ell}$-modules
 $$\Hom(X,Y)\otimes\Z_{\ell}\to \Hom(X(\ell),Y(\ell)).$$
 Explicitly, if $u \in \Hom(X,Y)\otimes\Z_{\ell}$ then for each
 $\nu$ we may pick
 $$w(\nu)\in \Hom(X,Y)=\Hom(X,Y)\otimes 1\subset
 \Hom(X,Y)\otimes\Z_{\ell}$$
 such that
 $$u-w(\nu) \in \ell^{\nu}\cdot
 \{\Hom(X,Y)\otimes\Z_{\ell}\}=\{\ell^{\nu}\cdot\Hom(X,Y)\}\otimes\Z_{\ell}=\Hom(X,Y)\otimes \ell^{\nu} \Z_{\ell}.$$
 Then the corresponding morphism of group schemes $u_{(\nu)}:=w(\nu): X_{\ell^{\nu}}
 \to Y$ does not depend on the choice of $w(\nu)$ and defines the
 corresponding morphism of $\ell$-divisible groups
 $$u_{(\nu)}: X_{\ell^{\nu}}\to Y_{\ell^{\nu}}; \ \nu=1, 2, \dots .$$

\begin{rem}
Since $\Hom(X,Y)$ is a free commutative group of finite rank, the
$\Z_{\ell}$-module
 $\Hom(X,Y)\otimes\Z_{\ell}$ is a free module of
 finite rank.
 \end{rem}

 The following assertion seems to be well known (at least, when
 $\ell\ne \fchar(K)$).

 \begin{lem}
The natural homomorphism of $\Z_{\ell}$-modules
$$\Hom(X,Y)\otimes\Z_{\ell}\to \Hom(X(\ell),Y(\ell))$$
is injective.
 \end{lem}

 \begin{proof}
If it is not injective and $u$ lies in the kernel then $u_{(\nu)} \in
\ell^{\nu}\cdot\Hom(X,Y)$ for all $\nu$. Since $u-u_{(\nu)}\in \ell^{\nu}\cdot
 \{\Hom(X,Y)\otimes\Z_{\ell}\}$, we conclude that
 $u\in \ell^{\nu}\cdot
 \{\Hom(X,Y)\otimes\Z_{\ell}\}$ for all $\nu$.  Since $\Hom(X,Y)\otimes\Z_{\ell}$ is a free $\Z_{\ell}$-module of
 finite rank, it follows that $u=0$.
 \end{proof}

 \begin{cor}
 \label{isoELL}
The following conditions are equivalent:

\begin{itemize}
\item[(i)]
 There exists an isogeny $u:X \to Y$, whose degree is not divisible by $\ell$.
 \item[(ii)]
 There exists $w \in \Hom(X,Y)\otimes\Z_{\ell}$ that induces an isomorphism of
 $\ell$-divisible group $X(\ell) \to Y(\ell)$.
\end{itemize}
 \end{cor}

\begin{proof}
Let $u:X \to Y$ be an isogeny, whose degree is not divisible by $\ell$.
Applying Lemma \ref{isogprime}(i)  to all $n=\ell^{\nu}$, we conclude that $u$
induces an isomorphism $X(\ell)\cong Y(\ell)$.

Now suppose that $w \in \Hom(X,Y)\otimes\Z_{\ell}$ that induces an isomorphism
of
 $\ell$-divisible group $X(\ell) \to Y(\ell)$. In particular, $w$ induces an
 isomorphism of finite group schemes $w_{(1)}:X_{\ell}\cong Y_{\ell}$. On the other
 hand, there exists $u \in \Hom(X,Y)$ such that
 $$w-u \in \ell \cdot \{\Hom(X,Y)\otimes\Z_{\ell}\}=\Hom(X,Y)\otimes\ell\Z_{\ell}.$$
 This implies that $u$ and $w$ induce the same morphism of finite group schemes
 $X_{\ell} \to Y_{\ell}$. It follows that the morphism
 $$u_{\ell}=u_{(1)}:X_{\ell} \to Y_{\ell}$$
 induced by $u$ coincides with $w_{(1)}$ and therefore is an isomorphism. Now Lemma
 \ref{isogprime}(ii) implies that $u$ is an isogeny, whose degree is not
 divisible by $\ell$.
\end{proof}

\end{sect}

\begin{sect}
\label{TateELL}
 {\bf Tate modules}. In this subsection we assume
that $\ell$ is a prime different from $\fchar(K)$. If $n=\ell^{\nu}$
then $X_n$ is an \'etale finite group scheme of order $n^{2\dim(X)}$
and we will identify its with the Galois module of its
$\bar{K}$-points. (Actually, all points of $X_n$ are defined over a
separable algebraic extension of $K$). The Tate $\ell$-module
$T_{\ell}(X)$ is defined as the projective limit of Galois modules
$X_{\ell^{\nu}}$ where the transition map $X_{\ell^{\nu+1}}\to
X_{\ell^{\nu}}$ is multiplication by $\ell$. The Tate module carries
a natural structure of free $\Z_{\ell}$-module of rank $2\dim(X)$;
it is also provided with a natural structure of Galois module in
such a way that natural homomorphisms $T_{\ell}(X)\to
X_{\ell^{\nu}}$ induce isomorphisms of Galois modules
$$T_{\ell}(X)\otimes \Z/\ell^{\nu} \cong X_{\ell^{\nu}}.$$
Explicitly, $T_{\ell}(X)$ is the set of all collections
$x=\{x_{\nu}\}_{\nu=1}^{\infty}$ with
$$x_{\nu}\in X_{\ell^{\nu}}, \quad x_{\nu+1}=\ell x_{\nu} \ \forall \nu .$$
The map $x\mapsto x_{\nu}$ defines the surjective homomorphism of
Galois modules $T_{\ell}(X)\to X_{\ell^{\nu}}$, whose kernel
coincides with $\ell^{\nu}\cdot T_{\ell}(X)$ and therefore induces
the isomorphism of Galois modules $T_{\ell}(X)/\ell^{\nu} \cong
X_{\ell^{\nu}}$ mentioned above.

 If $Y$ is an abelian variety over $K$ then we write
$\Hom_{\Gal}(T_{\ell}(X),T_{\ell}(Y))$ for the $\Z_{\ell}$-module
of all homomorphisms of $\Z_{\ell}$-modules $T_{\ell}(X) \to
T_{\ell}(Y)$ that commute with the Galois action(s), i.e., are
also homomorphisms of Galois modules.

The $\Z_{\ell}$-module $\Hom_{\Gal}(T_{\ell}(X),T_{\ell}(Y))$ is the set of
collections $w=\{w_{\nu}\}_{\nu=1}^{\infty}$ of homomorphisms of Galois modules
$$w_{\nu}:T_{\ell}(X)/\ell^{\nu}=X_{\ell^{\nu}}\to
Y_{\ell^{\nu}}=T_{\ell}(Y)/\ell^{\nu}$$  such that
$$w_{\nu}(x_{\nu})=\ell\cdot  uw_{\nu+1}(x_{\nu+1}) \quad \forall x=\{x_{\nu}\}_{\nu=1}^{\infty} \in
T_{\ell}(X).$$

Now if $z\in X_{\ell^{\nu}}$ then there exists $x \in T_{\ell}(X)$
with $x_{\nu}=z$. We have $\ell x_{\nu+1}=x_{\nu}=z$ and
$$w_{\nu}(z)=w_{\nu}(x_{\nu})=\ell\cdot
w_{\nu+1}(x_{\nu+1})=w_{\nu+1}(\ell
x_{\nu+1})=w_{\nu+1}(x_{\nu})=w_{\nu+1}(z),$$ i.e., the
restriction of $w_{\nu+1}$ to $X_{\ell^{\nu}}$ coincides with
$w_{\nu}$. This means that the collection
$\{w_{\nu}\}_{\nu=1}^{\infty}$ defines a morphism of
$\ell$-divisible groups over $K$
$$X(\ell) \to Y(\ell).$$

Conversely, if $u=\{u_{(\nu})\}_{\nu=1}^{\infty}$ is a morphism
$X(\ell) \to Y(\ell)$ over $K$ then $$u_{(\nu)}:X_{\ell^{\nu}}\to
Y_{\ell^{\nu}}$$ is a homomorphism of Galois modules; in addition,
the restriction of $u_{(\nu+1)}$ to $X_{\ell^{\nu}}$ coincides with
$u_{(\nu)}$. This implies that for each
$\{x_{\nu}\}_{\nu=1}^{\infty}\in T_{\ell}(X)$
$$u_{(\nu)}(x_{\nu})=u_{(\nu+1)}(x_{\nu})=u_{(\nu+1)}(\ell
x_{\nu+1})=\ell u_{(\nu+1)}(x_{\nu+1})$$ for all $\nu$. This means
that the collection $\{u_{(\nu)}\}_{\nu=1}^{\infty}$ defines a
homomorphism of Galois modules $T_{\ell}(X)\to T_{\ell}(Y)$. Those
observations give us the natural isomorphism of $\Z_{\ell}$-modules
$$\Hom(X(\ell),Y(\ell))=\Hom_{\Gal}(T_{\ell}(X),T_{\ell}(Y).$$
\end{sect}

\section{Useful results}
\label{use}

\begin{thm}[\cite{ZarhinMatZam,ZarhinInv85,MB}]
\label{principal} Let $X$ be an abelian variety of positive
dimension over a field $K$ and $X^t$ its dual. Then $(X\times
X^t)^4$ admits a principal $K$-polarization.
\end{thm}

We prove Theorem \ref{principal} in Section \ref{quatTRICK}.

\begin{thm}[\cite{LOZ}]
\label{subAV} Let $X$ be an abelian variety over $K$. The set of abelian
$K$-subvarieties of $X$ is finite, up to the action of the group $\Aut(X)$ of
$K$-automorphisms of $X$.
\end{thm}

 We sketch the proof of Theorem \ref{subAV} in Section \ref{pLOZ}.

\begin{lem}[Tate (\cite{Tate}, Sect. 2, p. 136)]
\label{poldegree}
 Let $K$ be a finite field, $g$ and $d$ are
positive integers. The set of $K$-isomorphism classes of $g$-dimensional
abelian varieties over $K$ that  admit a $K$-polarization of degree $d$ is
finite.
\end{lem}

Lemma \ref{poldegree} will be proven in Section \ref{polAV}.

\begin{thm}[\cite{ZarhinMatZam}, Th. 4.1]
\label{finitef} Let $K$ be a finite field, $g$ a positive integer.
Then the set of $K$-isomorphism classes of $g$-dimensional abelian
varieties over $K$ is finite.
\end{thm}

\begin{proof}[Proof of Theorem \ref{finitef} (modulo Theorem
\ref{principal} and Lemma \ref{poldegree})] Suppose that $X$ is a
$g$-dimensional abelian variety over $K$.
 By Lemma \ref{poldegree}, the set of
$8g$-dimensional abelian varieties over $K$ of the form $(X\times X^t)^4$ is
finite, up to $K$-isomorphism.  The abelian variety $X$ is isomorphic over $K$
to an abelian subvariety of $(X\times X^t)^4$. In order to finish the proof,
one has only to recall that thanks to Theorem \ref{subAV},  the set of abelian
subvarieties of a given abelian variety is finite, up to a $K$-isomorphism.
\end{proof}

 We need Theorem \ref{gabriel} in order to state the following
assertion.

\begin{cor}[Corollary to Theorem \ref{finitef}]
\label{semisimple}
 Let $X$ be  an abelian variety of positive
dimension over a finite field $K$. There exists a positive integer
 $r=r(X,K)$ that enjoys the following properties:

 \begin{itemize}
\item[(i)] If $Y$ is an abelian variety over $K$ that is
$K$-isogenous to $X$ then there exists a $K$-isogeny $\beta:X \to
Y$ such that $\ker(\beta)\subset X_r$.
 \item[(ii)]
 If $n$ is a positive integer and $W\subset X_n$ is a group
 subscheme over $K$ then there exists an endomorphism $u\in
 \End(X)$ such that
 $$r W \subset u X_n \subset W.$$
 \end{itemize}
\end{cor}

\begin{rem}
The assertion \ref{semisimple}(i) follows readily from Theorem
\ref{finitef}.
\end{rem}

We prove Corollary \ref{semisimple}(ii) in Section \ref{finiteKER}.

\section{Main results}
\label{mainres}


\begin{thm}
\label{endor} Let $X$ be an abelian variety of positive dimension
over a finite field $K$. There exists a positive integer
 $r_1=r_1(X,K)$ that enjoys the following properties:

 Let $n$ be a positive integer and $u_n \in \End(X_n)$. Let us put $m=n/(n,r_1)$. Then there
 exists $u\in \End(X)$  such that
 the images of $u$ and $u_n$ in $\End(X_m)$ do coincide.
 \end{thm}

We prove Theorem \ref{endor} in Section \ref{proofMAIN}.

 Applying Theorem \ref{endor} to a product $X=A\times B$ of
 abelian varieties $A$ and $B$, we obtain the following statement.

\begin{thm}
\label{homor} Let $A,B$ be  abelian varieties of positive
dimension over a finite field $K$. There exists a positive integer
 $r_2=r_2(A,B)$ that enjoys the following properties:

 Suppose that $n$ is a positive integer and $u_n:A_n \to B_n$ is a morphism
 of group schemes over $K$. Let us put $m=n/(n,r_2)$. Then there exists a
 homomorphism $u:A\to B$ of abelian varieties over $K$ such that
 the images of $u$ and $u_n$ in $\Hom(A_m,B_m)$ do coincide.
 \end{thm}

The following assertions follow readily from Theorem \ref{homor}.

\begin{cor}[First Corollary to Theorem \ref{homor}]
\label{primeL} If $n$ and $r_2$ are relatively prime (e.g., $n$ is
a prime that does not divide $r_2$) then the natural injection
$$\Hom(A,B)\otimes\Z/n \hookrightarrow \Hom(A_n,B_n)$$ is bijective.
\end{cor}

\begin{cor}[Second Corollary to Theorem \ref{homor}]
\label{primel} Let $\ell$ be a prime and $\ell^{r(\ell)}$ is the
exact power of $\ell$ dividing $r_2$. Then for each positive
integer $i$ the image of
$$\Hom(A_{\ell^{i+r(\ell)}},B_{\ell^{i+r(\ell)}})\to
\Hom(A_{\ell^{i}},B_{\ell^{i}})$$ coincides with the image of
$$\Hom(A,B)\otimes\Z/\ell^i \hookrightarrow \Hom(A_{\ell^i},B_{\ell^i}).$$
\end{cor}

\section{Abelian subvarieties}
\label{pLOZ}
 We follow the exposition in \cite{LOZ}.

 The next statement is a corollary of a finiteness result of
Borel and Harish-Chandra \cite[Theorem 6.9]{Borel}; it  may  also be
deduced from the Jordan--Zassenhaus theorem \cite[Theorem
26.4]{Reiner}.

\begin{prop}[\cite{LOZ}, p. 514]
\label{LOZf}
 Let $F$ be a finite-dimensional semisimple $\Q$-algebra, $M$ a
finitely generated right $F$-module, $L$ a $\Z$-lattice in $M$. Let $G$ be the
group of those automorphisms $\sigma$ of the $F$-module $M$ for which
$\sigma(L)=L$. Then the number of $G$-orbits of the set of $F$-submodules of
$M$ is finite.
\end{prop}

Now let $X$ be an abelian variety over $K$. We are going to apply Proposition
\ref{LOZf} to
$$F=\End(X)\otimes\Q,\ M=\End(X)\otimes\Q, \ L=\End(X).$$
One may identify $G$ with the group $\Aut(X)=\End(X)^{*}$ of
automorphisms of $X$: here elements of $\End(X)^{*}$ act as left
multiplications on $\End(X)\otimes\Q=M$.

On the other hand, to each abelian $K$-subvariety $Y\subset X$ corresponds the
right ideal $$I(Y)=\{u \in \End(X)\mid u(X)\subset Y\}$$ and the $F$-submodule
$$I(Y)_{\Q}=I(Y)\otimes\Q \subset \End(X)\otimes\Q=M.$$
Using the theorem of Poincar\'e--Weil \cite[Proposition
12.1]{Milne}, one may prove (\cite[p. 515]{LOZ} that $I(Y)_{\Q}$
uniquely determines $Y$. Even better, if $Y'$ is an abelian
$K$-subvariety of $X$ and
$$ u I(Y)_{\Q}=I(Y')_{\Q}$$
for $u \in \Aut(X)=\End(X)^{*}$ then $Y'=u(Y)$. Now Proposition \ref{LOZf}
implies the finiteness of the number of orbits of the set of abelian
$K$-subvarieties of $X$ under the natural action of $\Aut(X)$. This proves
Theorem \ref{subAV}. (See  \cite{LOZ2} for variants and complements.)

\section{Polarized abelian varieties}
\label{polAV}

\begin{lem}[Mumford's lemma \cite{GIT}]
\label{MumL} Let $X$ be an abelian variety of positive dimension over a field
$K$. If $\lambda: X \to X^t$ is a polarization then there exists an ample
invertible sheaf $\L$ on $X$ such that
$$\Lambda_{\bar{\L}}=2\bar{\lambda}$$
where $\bar{\L}$ is the  invertible sheaf on $\bar{X}$ induced by
${\L}$.
\end{lem}

\begin{proof}
See \cite[Ch. 6, Sect. 2, pp. 120--121]{GIT} where a much more
general case of abelian schemes is considered. (In notation of
\cite{GIT}, $S$ is the spectrum of $K$.) Let me just recall an
explicit construction of $\L$. Let $\P$ be the universal
Poincar\'e invertible sheaf on $X\times X^t$ \cite[Sect.
9]{Milne}. Then $\L:=(1_X,\lambda)^{*}\P$ where $(1_X,\lambda):X
\to X\times X^t$ is defined by the formula
$$x \mapsto (x,\lambda(x)).$$
\end{proof}

\begin{proof}[Proof of Lemma \ref{poldegree}]
So, let $X$ be a $g$-dimensional abelian variety over a finite field $K$ and
let $\lambda: X \to X^t$ be a polarization of degree $d$.
 We follow the
exposition in \cite[p. 243]{Ram}. By Lemma \ref{MumL}, there exists
an invertible ample sheaf $\L$ on $X$ such that the
self-intersection index of $\bar{\L}$ equals $2^g d g!$ \cite[Sect.
16]{MumfordAV}. The invertible sheaf $\bar{\L}^3$ is very ample, its
space of global section has dimension $6^g d$; the self-intersection
index of $\L$ equals $6^g d g!$ \cite[Sect. 16]{MumfordAV}. This
implies that $\L^3$ is also very ample and gives us an embedding
(over $K$) of $X$ into the $6^g d-1$-dimensional projective space as
a closed $K$-subvariety of degree $6^g d g!$. All those subvarieties
are uniquely determined by their Chow forms (\cite[Ch. 1, Sect.
6.5]{Sha}, \cite[Lecture 21, pp. 268--273]{H}), whose coefficients
are elements of $K$. Since $K$ is finite and the number of
coefficients depends only on the degree and dimension, we get the
desired finiteness result.
\end{proof}

\section{Quaternion trick}
\label{quatTRICK}
 Let $X$ be an abelian variety of positive dimension over a field $K$.
and $\lambda: X \to X^t$  a $K$-polarization. Pick a positive integer $n$ such
that
$$\ker(\lambda)\subset X_n.$$

\begin{lem}
Suppose that there exists an integer $a$ such that $a^2+1$ is divisible by $n$.
Then $X\times X^t$ admits a principal polarization that is defined over $K$.
\end{lem}

\begin{proof}
Let $$V\subset \ker(\lambda)\times \ker(\lambda)\subset X_n\times
X_n\subset X \times X$$ be the graph of multiplication by $a$ in
$\ker(\lambda)$. Clearly, $V$ is a finite group subscheme over $K$
that is isomorphic to $\ker(\lambda)$  and therefore its order is
equal to $\deg(\lambda)$. Notice that $\deg(\lambda)$ is the
square root of $\deg(\lambda^2)$.

For each commutative $\bar{K}$-algebra $R$ the group $\bar{V}(R)$ of
$R$-points coincides with the set of all the pairs $(x,ax)$ with
$x\in \ker(\bar{\lambda})\subset \bar{X}_n$. This implies that for
all $(x,ax),(y,ay)\in \bar{V}(R)$ we have
$$e_{\lambda^2}((x,ax),(y,ay))=e_{\lambda}(x,y)\cdot e_{\lambda}(ax,ay)=
e_{\lambda}(x,y)\cdot e_{\lambda}(a^2 x,y)=$$
$$e_{\lambda}(x,y)\cdot
e_{\lambda}(-x,y)=e_{\lambda}(x,y)/e_{\lambda}(x,y)=1.$$ In other
words, $\bar{V}$ is isotropic with respect to $e_{\lambda^2}$; in
addition, $$\#(\bar{V})^2=\deg(\lambda)^2=\deg(\lambda^2).$$ This
implies that $X^2/V$ is a principally polarized abelian variety
over $K$. On the other hand, we have an isomorphism of abelian
varieties over $K$
$$f: X \times X \to X \times X=X^2, \ (x,y)\mapsto
(x,ax)+(0,y)=(x,ax+y)$$ and
$$V=f(\ker{\lambda}\times\{0\})\subset f(X\times \{0\}).$$
Thus, we obtain $K$-isomorphisms
$$X^2/V \cong X/\ker(\lambda)\times X=X^t\times X=X\times X^t.$$
In particular, $X\times X^t$ admits a principal $K$-polarization
and we are done.
\end{proof}

\begin{proof}[Proof of Theorem \ref{principal}]
Choose a quadruple of integers $a,b,c,d$ such that $$0\ne
s:=a^2+b^2+c^2+d^2$$ is congruent to $-1$ modulo $n$. We denote by
$\I$ the ``quaternion"
$$\I=
\begin{pmatrix}
a & -b & -c & -d\\
b & \ a & \ d &\  c\\
c & -d & \ a & \ b\\
d & \ c & -b & \ a\\
\end{pmatrix}
\in \Mat_4(\Z)\subset \Mat_4(\End(X)=\End(X^4).$$ We have
$$\I^{*} \I=a^2+b^2+c^2+d^2=s \in \Z \subset \Mat_4(\Z)\subset \Mat_4(\End(X)=\End(X^4).$$
Let $$V\subset \ker(\lambda^4)\times \ker(\lambda^4)\subset
(X^4)_n\times (X^4)_n\subset X^4 \times X^4=X^8$$ be the graph of
$$\I: \ker(\lambda^4) \to \ker(\lambda^4).$$
 Clearly, $V$ is a finite
group subscheme over $K$ and its order is equal to
$\deg(\lambda^4)$. Notice that $\deg(\lambda^4)$ is the square
root of $\deg(\lambda^8)$.

For each commutative $\bar{K}$-algebra $R$ the group $\bar{V}(R)$ of
$R$-points consists of all the pairs $(x,\I x)$ with $x\in
\ker(\bar{\lambda}^4)\subset (\bar{X^4})_n$. This implies that for
all $(x,\I x),(y,\I y)\in \bar{V}(R)$ we have
$$e_{\lambda^4}((x,\I x),(y,\I y))=e_{\lambda^4}(x,y)\cdot e_{\lambda^4}(\I x,\I y)=
e_{\lambda^4}(x,y)\cdot e_{\lambda}(x,\I^t \I y)=$$
$$e_{\lambda}(x,y)\cdot
e_{\lambda}(x,s y)=e_{\lambda}(x,y)\cdot e_{\lambda}(x,-y)
=e_{\lambda}(x,y)/e_{\lambda}(x,y)=1.$$ In other words, $\bar{V}$
is isotropic with respect to $e_{\lambda^4}$; in addition,
$$\#(\bar{V})^2=\deg(\lambda^4)^2=\deg(\lambda^8).$$ This implies that
$X^8/V$ is a principally polarized abelian variety over $K$. On
the other hand, we have an isomorphism of abelian varieties  over
$K$
$$f: X^4 \times X^4 \to X^4 \times X^4=X^8, \ (x,y)\mapsto
(x,\I x)+(0,y)=(x,\I x+y)$$ and
$$V=f(\ker(\lambda^4)\times\{0\})\subset f(X^4\times \{0\}).$$
Thus, we obtain $K$-isomorphisms
$$X^4/V \cong X^4/\ker{\lambda^4}\times X^4=(X^4)^t\times X^4=(X\times X^t)^4.$$
In particular, $(X\times X^t)^4$ admits a principal
$K$-polarization and we are done.
\end{proof}

\begin{rem}
We followed the exposition in \cite[Lemma 2.5]{ZarhinMatZam},
\cite[Sect. 5]{ZarhinInv85}. See \cite[Ch. IX, Sect. 1]{MB} where
Deligne's proof is given.
\end{rem}

\section{Finite group subschemes of abelian varieties}
\label{finiteKER}
\begin{proof}[Proof of Corollary \ref{semisimple}(ii)]
Let $r$ be as in \ref{semisimple}(i).  Let us consider the abelian
variety $Y:=X/W$ and the canonical isogeny $K$-isogeny  $\pi:X \to
X/W=Y$. Clearly,
$$W=\ker(\pi).$$
Since $W\subset X_n$, there exists a $K$-isogeny $v:Y\to X/X_n=X$ such that the
composition $v\pi$ coincides with multiplication by $n$ in $X$; in addition,
$$\pi n_X=n_Y\pi:X \to Y$$
is a $K$-isogeny, whose degree is  $\#(W)\times n^{2\dim(X)}$.
Here $n_X$ (resp. $n_Y$) stands for multiplication by $n$ in $X$
(resp. in $Y$).
                Let us put
$$U=\ker(\pi n_X)=\ker(n_Y\pi)\subset X;$$
it is a finite commutative group $K$-(sub)scheme and
$$\#(U)=\#(W)\times n^{2\dim(X)}.$$
Then
$$X_n\subset U, \ W \subset U; \ \pi(U)\subset Y_n, \ n_X(U)\subset W.$$
The order arguments imply that the natural morphisms of group
$K$-schemes $$\pi: U \to Y_n, \ n_X: U \to W$$ are surjective,
i.e.,
$$\pi(U)= Y_n, \ n U = W.$$
We have
$$v(Y_n)=v(\pi(U))=v\pi (U)=n U=W,$$
i.e.,
$$v(Y_n)=W.$$

By \ref{semisimple}(i), there exists a $K$-isogeny $\beta: X\to Y$
with $\ker(\beta)\subset X_r$.  Then there exists a $K$-isogeny
$\gamma:Y\to X$ such that $\gamma\beta=r_X$. This implies that
$$\gamma r_Y=r_X \gamma=\gamma\beta \gamma=\gamma(\beta \gamma),$$
i.e.,
$$\gamma r_Y=\gamma(\beta \gamma).$$
It follows that $r_Y=\beta \gamma$, because $\ker(\gamma)$ is
finite while $(r_Y-\beta \gamma)Y$ is an abelian subvariety. This
implies that
$$\beta (X_n) \supset \beta(\gamma(Y_n))= \beta \gamma(Y_n)=r Y_n.$$
Let us put
$$u=v\beta\in \End(X).$$
We have
$$Y_n\supset \beta (X_n)\supset r Y_n.$$
This implies that
$$W=v(Y_n)\supset v(\beta)(X_n)=u(X_n),$$
$$u(X_n)=v(\beta(X_n))\supset v(r Y_n)=r(W)$$
and therefore
$$W\supset u(X_n)\supset r (W).$$
\end{proof}

\section{Dividing homomorphisms of abelian varieties}
\label{finiteDIV} Results of this Section will be used in the proof of Theorem
\ref{endor} in Section \ref{proofMAIN}.

 Throughout this Section, $Y$ is an abelian variety over a field $K$.
 The following statement is  well known.

\begin{lem}
\label{divprin} let $u: Y\to Y$ be a $K$-isogeny. Suppose that $Z$ is an
abelian variety over $K$. Let $v \in \Hom(Y,Z)$ and $\ker(u)\subset \ker(v)$
(as a group subscheme in $Y$). Then there exists exactly one $w \in \Hom(Y,T)$
such that $v=wu$, i.e., the diagram

\begin{displaymath}
\xymatrix{
   Y \ar[r]^{u} \ar[dr]_{v} & Y \ar[d]^{w} \\
                                       & Z
}
\end{displaymath}
is commutative. In addition, $w$ is an isogeny if and only if $v$
is an isogeny.
\end{lem}

\begin{proof}
We have $Y\cong Y/\ker(u)$. Now the result follows from the
universality property of quotient maps.

\end{proof}

Let $n$ be  a positive integer and $u$ an endomorphism of $Y$. Let
us consider the homomorphism of abelian varieties over $K$
$$(n_Y,u): Y \to Y\times Y, \quad y\mapsto (ny,uy).$$
Then
$$\ker((n_Y,u))=\ker(Y_n \stackrel{u}{\to} Y_n)\subset Y_n\subset Y.$$
Slightly abusing notation, we denote the finite commutative group
$K$-(sub)scheme $\ker((n_Y,u))$ by $\{\ker(u)\bigcap Y_n\}$.

\begin{lem}
\label{div} Let $Y$ be an abelian variety of positive dimension
over a  field $K$. Then there exists a positive integer $h=h(Y,K)$
that enjoys the following properties:

If $n$ is  a positive integer, $u,v \in \End(Y)$  are
endomorphisms such that
$$\{\ker(u)\bigcap Y_n\}\subset \{\ker(v)\bigcap Y_n\}$$
then there exists a $K$-isogeny $w: Y\to Y$ such that
$$hv-wu \in n\cdot \End(Y).$$
In particular, the images of $hv$ and $wu$ in  $\End(Y_n)$ do
coincide.
\end{lem}

\begin{proof}
Since $\OC:=\End(Y)$ is an order in the semisimple
finite-dimensional $\Q$-algebra $\End(Y)\otimes\Q$, the
Jordan--Zassenhaus  theorem \cite[Th. 26.4]{Reiner} implies that
there exists a positive integer $M$ that enjoys the following
properties:

{\sl if $I$ is a left ideal in $\OC$ that is also a subgroup of finite index
then there exists $a_I \in \OC$ such that the principal left ideal $a\cdot \OC$
is a subgroup in $I$ of finite index dividing $M$; in particular},
$$ M \cdot I\subset a_I \cdot \OC \subset I.$$
Clearly, such $a_I$ is invertible in $\End(Y)\otimes\Q$ and
therefore is an isogeny. Let us put
$$h:=M^3.$$

Let us consider the left ideals
$$I=n\OC+u\OC, \ J=n\OC+v\OC$$
in  $\OC$. Then both $I$ and $J$ are subgroups of finite index in $\OC$. So,
there exist $K$-isogenies
$$a_I: Y\to Y, \ a_J: Y\to Y$$ such that
$$M \cdot I\subset a_I \cdot \OC \subset I, \ M \cdot I\subset a_J \cdot \OC \subset
J.$$ In particular, there exist $b,c \in \OC$ such that
$$Ma_I-bu  \in n\cdot \OC, \ Mv=c a_J.$$
In obvious notation
$$\{\ker(v)\bigcap Y_n \} \subset    \ker(a_J)\subset \{\ker(Mv)\bigcap Y_{Mn}\}=M^{-1}\{\ker(v)\bigcap
Y_n\}\subset Y,$$
$$\{\ker(u)\bigcap Y_n\} \subset    \ker(a_I)\subset \{\ker(Mu)\bigcap Y_{Mn}\}=M^{-1}\{\ker(u)\bigcap
Y_n\}\subset Y.$$ This implies that
$$\ker(a_I)\subset M^{-1}\{\ker(u)\bigcap Y_n\}\subset M^{-1}\{\ker(v)\bigcap
Y_n\}\subset M^{-1}\ker{(a_J})=\ker(M a_J)$$ and therefore
$$\ker(a_I)\subset \ker(M a_J).$$
By Lemma \ref{divprin}, there exists a $K$-isogeny $z: Y \to Y$
such that $M a_J =z a_I$ and therefore $M^2 a_J =M z a_I$. This
implies that
$$M^3 v=M^2 c a_J=M c (M a_J)=M c (z a_I)= c z (M a_I)=$$
$$cz [bu+(M a_I -bu)]=(czb)u+ cz(M a_I -bu).$$
Since $h=M^3$ and $bu - M a_I$ is divisible by $n$ in
$\OC=\End(Y)$,
$$h v -(czb)u \in n\cdot \End(Y).$$
So, we may put  $w=czb$.
\end{proof}

\section{Endomorphisms of group schemes}
\label{proofMAIN}
\begin{proof}[Proof of Theorem \ref{endor}]
Let $X$ be an abelian variety of positive dimension over a finite
field $K$. Let us put $Y:=X\times X$. Let $h=h(Y)$ be as in Lemma
\ref{div} and $r=r(Y,K)$ be as in Corollary \ref{semisimple}. Let
us put
$$r_1=r_1(X,K):= r(Y,K) h(Y,K).$$
Let $n$ be a positive integer and $u_n\in \End(X_n)$. Let $W$ be
the graph of $u_n$ in $X_n\times X_n=(X\times X)_n=Y_n$, i.e., the
image of
$$({\mathbf 1}_n,u_n): X_n \hookrightarrow X_n\times X_n=(X\times
X)_n=Y_n.$$ Here ${\mathbf 1}_n$ is the identity automorphism of
$X_n$.

By Corollary \ref{semisimple}, there exists $v \in \End(Y)$ such
that
$$r W \subset u(Y_n)\subset W.$$

 Let $\pr_1, \pr_2: Y=X\times X \to X$ be the projection maps and
$$q_1: X=X\times \{0\}\subset X\times X=Y, \ q_2: X=\{0\}\times X\subset X\times
X=Y$$ be the inclusion maps. Let us consider the
homomorphisms
$$\pr_1 v, \pr_2 v : Y \to X$$ and the endomorphisms
$$v_1=q_1 \pr_1 v, \ v_2=q_1 \pr_2 v \in \End(X\times
X)=\End(Y).$$ Clearly, $$v:Y\to Y=X\times X$$ is ``defined" by
pair
$$(\pr_1 v, \pr_2 v): Y \to X\times X=Y.$$
Since $W$ is a graph,
$$\pr_1(W)=X_n, \  v(Y_n)\subset W$$
and
$$\{\ker(\pr_1 v)\bigcap Y_n\} \subset \{\ker(\pr_2 v)\bigcap Y_n\}.$$
Since $q_1$ and $q_2$ are embeddings,
$$\{\ker(v_1)\bigcap Y_n\} \subset \{\ker(v_2)\bigcap Y_n\}.$$
 By Lemma \ref{div}, there exists a $K$-isogeny $w:Y\to Y$
such that the restrictions of $hv_2$ and $w v_1$ to $Y_n$ do
coincide. Taking into account that
$$v_1(X\times X)\subset X\times \{0\}, \ v_2(X\times X)\subset \{0\}\times
X,$$ we conclude that if we put
$$w_{12}=\pr_2 w q_1\in \End(X)$$
then the images of $h \ \pr_2 v$ and  $w_{12} \pr_1 v$ in
$\Hom(Y_n,X_n)=\Hom(X_n\times X_n,X_n)$ do coincide.

Since $W$ is the graph of $u_n$ and $u(Y_n)\subset W$,
$$\pr_2 v=u_n \pr_1 v \in \Hom(Y_n,X_n);$$
here both sides are viewed as  morphisms of group schemes $Y_n \to
X_n$. This implies that in $\Hom(Y_n,X_n)$ we have
$$w_{12} \pr_1 v = h \ \pr_2 v= h \ u_n \pr_1 v.$$
This implies that $w_{12}=h \ u_n$ on
$$\pr_1 v(Y_n)\subset X_n.$$
We have
$$\pr_1 v(Y_n)\supset r\  \pr_1(r(W))=r(X_n)$$
and therefore $w_{12}=h \ u_n$ on $r(X_n)$. By Lemma \ref{Xn},
$$r(X_n)=X_{n_1},$$
where $n_1=n/(n,r)$. So, $w_{12}=h\ u_n$ on $X_{n_1}$. Let us put
$d:=(n_1,h)$. Clearly, $X_d \subset X_{n_1}$ and $w_{12}=h u_n$
kills $X_d$, because $d$ divides $h$. This implies that there
exists $u \in \End(X)$ such that $w_{12}=d \ u$.  If we put
$m=n_1/d$ then  $h/d$ is a positive integer relatively prime to
$m$ and $(h/d) \ u \ d = (h/d) \ u_n\ d$ on $X_{n_1}$ and
therefore $(h/d)\ u = (h/d)\ u_n$ on $d(X_{n_1})=X_m$. Since
multiplication by $(h/d)$ is an automorphism of $X_m$, we conclude
that $u=u_n$ on $X_m$.

\end{proof}

\begin{cor}
\label{infS} Let $K$ be a finite field,  $X$ and $Y$ abelian varieties over
$K$. Let $S$ be the set of positive integers $n$ such that the finite
commutative group $K$-schemes $X_n$ and $Y_n$ are isomorphic. If $S$ is
infinite then $X$ and $Y$ are isogenous over $K$. In addition, if $S$ is the
set of powers of a prime $\ell$ then there exists a $K$-isogeny $X\to Y$, whose
degree is not divisible by $\ell$.
\end{cor}

\begin{proof}
Pick $n \in S$ such that $n>r_2:=r_2(X,Y)$ where $r_2$ is as in Theorem
\ref{homor}. Then $m:=n/(n,r_2)$ is strictly greater than $1$. (In addition, if
$n$ is a power of $\ell$ then $m$ is also a power of $\ell$.( Fix an
isomorphism $w_n: X_n \cong Y_n$. By Theorem \ref{homor}, there exists $u\in
\Hom(X,Y)$ such that the induced morphism $u_m:X_m \to Y_m$ coincides with the
restriction (image) of $w_n$ to (in) $\Hom(X_m,Y_m)$. But this restriction is
an isomorphism, since $w_n$ is an isomorphism. It follows that $u_m$ is an
isomorphism. Now the desired result follows from Lemma \ref{isogprime}(ii).

\end{proof}

\begin{thm}[Tate's theorem on homomorphisms]
 \label{tatehomo} Let $K$ be a finite field, $\ell$ an arbitrary
prime, $X$ and $Y$ abelian varieties over $K$ of positive
dimension. Let $X(\ell)$ and $Y(\ell)$ be the $\ell$-divisible
groups attached to $X$ and $Y$ respectively. Then the natural
embedding
$$\Hom(X,Y)\otimes\Z_{\ell} \hookrightarrow
\Hom(X(\ell),Y(\ell))$$ is bijective.
\end{thm}

\begin{rem}
Our proof will work for both cases $\ell\ne\fchar(K)$ and
$\ell=\fchar(K)$.
\end{rem}

\begin{proof}[Proof of Theorem \ref{tatehomo}]
Any element of $\Hom(X(\ell),Y(\ell))$ is a collection $$\{w_{(\nu)}
\in \Hom(X_{\ell^{\nu}},Y_{\ell^{\nu}})\}_{{\nu}=1}^{\infty}$$ such
that every $w_{(\nu)}$ coincides with the ``restriction" of
$w_{({\nu}+1)}$ to $X_{\ell^{\nu}}$. It follows from Corollary
\ref{primel} that there exists $u_{\nu}\in \Hom(X,Y)\otimes
\Z/\ell^{\nu}$ such that $w_{(\nu)}=u_{\nu}$. This implies that the
image of $u_{{\nu}+1}$ in $\Hom(X,Y)\otimes \Z/\ell^{{\nu}}$
coincides with $u_{\nu}$ for all ${\nu}$. This means that if $u$ is
the projective limit of $u_{\nu}$ in $\Hom(X,Y)\otimes \Z_{\ell}$
then $u$ induces (for all ${\nu}$) the morphism from
$X_{\ell^{\nu}}$ to $Y_{\ell^{\nu}}$ that coincides with $u_{\nu}$
and therefore with $w_{(\nu)}$.
\end{proof}

\begin{cor}
\label{isomDIV} Let $K$ be a finite field, $\ell$ an arbitrary prime, $X$ and
$Y$ abelian varieties over $K$ of positive dimension. Then the following
conditions are equivalent:

\begin{itemize}
\item There exists a $K$-isogeny $X \to Y$, whose degree is not divisible by
$\ell$.

\item The $\ell$-divisible groups $X(\ell)$ and $Y(\ell)$ are isomorphic.
\end{itemize}
\end{cor}

\begin{proof} It follows readily from Theorem \ref{tatehomo} and Corollary
\ref{isoELL}.
\end{proof}

\section{Homomorphisms of Tate modules and isogenies}
\label{TateMOD}
 Throughout this Section, $K$ is a finite field and $\ell$ is a prime
$\ne \fchar(K)$.

 Combining  Theorem \ref{tatehomo} with results of Section
\ref{TateELL}, we obtain the following statement.

\begin{thm}[Tate \cite{Tate}]
\label{tateTl}  Let $X$ and $Y$ be abelian varieties over $K$. Then
$$\Hom(X,Y)\otimes
\Z_{\ell}= \Hom_{\Gal}(T_{\ell}(X),T_{\ell}(Y)).$$
\end{thm}

Let $X$ be an abelian variety over $K$. Let us consider the $\Q_{\ell}$-vector
space
$$V_{\ell}(X)=T_{\ell}(X)\otimes_{\Z_{\ell}}\Q_{\ell}$$
provided with the natural structure of Galois module. We have
$$\dim_{\Q_{\ell}}(V_{\ell}(X))=2\dim(X)$$
and the map
$$T_{\ell}(X)\hookrightarrow V_{\ell}(X), \ z \mapsto z\otimes 1$$
identifies $T_{\ell}(X)$ with a Galois-invariant $\Z_{\ell}$-lattice. This
implies that  the natural map
$$\Hom_{\Gal}(T_{\ell}(X),T_{\ell}(Y))\otimes_{\Z_{\ell}}\Q_{\ell}\to
  \Hom_{\Gal}(V_{\ell}(X),V_{\ell}(Y))$$
  is bijective. Here $\Hom_{\Gal}(V_{\ell}(X),V_{\ell}(Y))$ is the
  $\Q_{\ell}$-vector space of $\Q_{\ell}$-linear homomorphisms of Galois
  modules $V_{\ell}(X)\to V_{\ell}(Y)$.

  Applying Theorem \ref{tateTl}, we obtain the following
  statement.

\begin{thm}[Tate \cite{Tate}]
\label{tateVl}  Let $X$ and $Y$ be abelian varieties over $K$. Then the natural
map
$$\Hom(X,Y)\otimes
\Q_{\ell}= \Hom_{\Gal}(V_{\ell}(X),V_{\ell}(Y))$$ is bijective.
\end{thm}

The following assertion is very useful.

\begin{cor}
[Tate's isogeny theorem \cite{Tate}] \label{tateIl}  Let $X$ and $Y$ be abelian
varieties over $K$. Then $X$ and $Y$ are isogenous over $K$ if and only if the
Galois modules $V_{\ell}(X)$ and $V_{\ell}(Y)$ are isomorphic.
\end{cor}

\begin{proof}
If $X$ and $Y$ are isogenous over $K$ then there exist a positive integer $N$
and isogenies $$\alpha: X\to Y, \ \beta:Y\to X$$ such that
$$\beta\alpha=N_X, \ \alpha\beta =N_Y.$$
By functoriality, $\alpha$ and $\beta$ induce homomorphisms of Galois modules
$$\alpha{(\ell)}:V_{\ell}(X)\to V_{\ell}(Y), \ \beta{(\ell)}:V_{\ell}(Y)\to
V_{\ell}(X)$$ such that the compositions
$\beta{(\ell)}\alpha{(\ell)}$ and $\alpha{(\ell)}\beta{(\ell)}$
coincide with multiplication by $N$ in $V_{\ell}(X)$ and
$V_{\ell}(Y)$ respectively. It follows that  $\alpha{(\ell)}$ is an
isomorphism of Galois modules $V_{\ell}(X)$ and $V_{\ell}(Y)$.

Suppose now that the Galois modules $V_{\ell}(X)$ and $V_{\ell}(Y)$ are
isomorphic. Then their $\Q_{\ell}$-dimensions coincide and therefore
$$\dim(X)=\dim(Y).$$
Choose an isomorphism
$$w:V_{\ell}(X)\cong V_{\ell}(Y)$$
of Galois modules. Replacing (if necessary) $w$ by $\ell^M w$ for
sufficiently large positive integer $M$, we may and will assume that
$$w(T_{\ell}(X))\subset T_{\ell}(Y).$$
The image $w(T_{\ell}(X))$ is a $\Z_{\ell}$-lattice in
$V_{\ell}(Y)$. This implies that $w(T_{\ell}(X))$ is a subgroup of
finite index in $T_{\ell}(Y)$. So, we may view $w$ as an {\sl
injective} homomorphism $T_{\ell}(X)\to T_{\ell}(Y)$ of Galois
modules. There exists a positive integer $M$ such that if
$$w'\in \Hom_{\Gal}(T_{\ell}(X),T_{\ell}(Y)), \ w'-w \in \ell^M \cdot
\Hom_{\Gal}(T_{\ell}(X),T_{\ell}(Y))$$ then
$$w':T_{\ell}(X)\to  T_{\ell}(Y)$$ is also injective.
 Since $\Hom(X,Y)$ is everywhere
dense with respect to $\ell$-adic topology in
$$\Hom(X,Y)\otimes
\Z_{\ell}= \Hom_{\Gal}(T_{\ell}(X),T_{\ell}(Y)),$$ there exists $u \in
\Hom(X,Y)$ such that the induced (by $u$) homomorphism of Galois modules
  $$u{(\ell)}:T_{\ell}(X)\to T_{\ell}(Y)$$
  is injective. This implies that
$$\rk_{\Z_{\ell}}(u{(\ell)}(T_{\ell}(X)))=\rk_{\Z_{\ell}}(T_{\ell}(X))=2\dim(X)=2\dim(Y).$$
   I claim that $u$ is an isogeny. Indeed, let us put
  $Z:=u(X)$: it is a (closed)  abelian subvariety of $Y$ that is defined over $K$.
  The homomorphism $u:X \to Y$ coincides with the composition of the natural surjection
  $X \to Z$ and the inclusion map $j:Z \hookrightarrow X$. This implies that
  $u{(\ell)}(T_{\ell}(X))$ is contained in $j{(\ell)}(T_{\ell}(Z))$
  where
  $$j(\ell):T_{\ell}(Z) \to T_{\ell}(Y)$$
is the homomorphism of Tate modules induced by $j$.
 It follows
  that
  $$2\dim(Z)=\rk(T_{\ell}(Z))\ge
  \rk(j(\ell)(T_{\ell}(Z)))\ge$$
  $$\rk(u{(\ell)}(T_{\ell}(X)))=2\dim(X)=2\dim(Y)$$
  and therefore $\dim(Z)\ge \dim(Y)$. (Hereafter $\rk$ stands for the rank of a
  free $\Z_{\ell}$-module.)

   Since $Z$ is a closed subvariety of $Y$,
  we conclude that $\dim(Z)=\dim(Y)$ and therefore $Z=Y$. This implies that
  $u: X \to Y$ is surjective. Since $\dim(X)=\dim(Y)$, we conclude that $u$ is
  an isogeny.
\end{proof}

Corollary \ref{tateIl} admits the following  ``refinement".

\begin{cor}
Let $X$ and $Y$ be abelian varieties over $K$. The following assertions are
equivalent.

\begin{itemize}
\item There exists an isogeny $X \to Y$, whose degree is not divisible by
$\ell$.

\item
 The Galois
modules $T_{\ell}(X)$ and $T_{\ell}(Y)$ are isomorphic.
\end{itemize}
\end{cor}

\begin{proof}
It follows readily from Corollary \ref{isomDIV} and the last displayed formula
in Subsection \ref{TateELL}.
\end{proof}

\section{An example}
\label{nonISOM}

Corollaries \ref{infS} and Corollary \ref{isomDIV}  suggest the following
question: if $X$ and $Y$ are abelian varieties over a finite field $K$ such
that $X_n \cong Y_n$ for all $n$ and $X(\ell) \cong Y(\ell)$ for all $\ell$
then is it true that $X$ and $Y$ are isomorphic? The aim of this Section is to
give a negative answer to this question. Our construction is based on the
theory of elliptic curves with complex multiplication \cite{SerreCM,LangE}.

We start to work over the field $\C$ of complex numbers.
 Let $F\subset \C$ be an imaginary quadratic field with the ring of integers
$\OC_F$. For every non-zero ideal $\b\subset \OC_F$ there exists an elliptic
curve $E^{(\b)}$ over $\C$ such that that its group of complex points
$E^{(\b)}(\C)$ (viewed as a complex Lie group) is $\C/\b$. There is a natural
ring isomorphism $\OC_F \cong \End(E^{(\b)})$ where any $a \in \OC_F$ acts on
$E^{(\b)}(\C)$ as
$$z+\b \mapsto az+\b.$$
In particular, $E^{(\b)}$ is an elliptic curve with complex
multiplication and $\j(E^{(\b)})\in \C$ is an {\sl algebraic
integer}.

Let us put $E:=E^{(\OC_F)}$. There is a natural bijection of groups
$$\b \cong \Hom(E,E^{(\b)}), \ c \mapsto u(c),$$
where homomorphism $u(c)$ acts on complex points as
$$u(c): \C/\OC_F \to \C/\b, \ z+\OC_F \mapsto c z+ \b.$$
In addition, for every non-zero $c$ the homomorphism $u(c):E\to E^{(\b)}$ is an
isogeny, whose  degree is the order of the (finite) quotient $\b/c \OC_F$. In
particular, $E$ and $E^{(\b)}$ are isomorphic if and only if $\b$ is a
principal ideal. This implies that if $\b$ is {\sl} not principal then
$$\j(E^{(\b)}) \ne \j(E).$$

\begin{lem}
For every prime $\ell$ there exists a non-zero $c \in \b$ such that the order
of $\b/c \OC_F$ is not divisible by $\ell$.
\end{lem}

\begin{proof}
We may assume that $\b$ is not principal.  If $\ell \OC_F$ is a prime ideal in
$\OC_F$, pick any $c \in \b \setminus \ell\b$. If $\ell \OC_F$ is a square
${\mathfrak L}^2$ of a prime ideal ${\mathfrak L}$, pick any $c \in
\b\setminus{\mathfrak L}\cdot \b$.
 If $\ell \OC_F$ is a product ${\mathfrak L}_1  {\mathfrak
L}_2$ of two distinct prime ideals ${\mathfrak L}_1, {\mathfrak L}_2
\subset \OC_F$, pick
$$c_1 \in {\mathfrak L}_1\cdot b\setminus {\mathfrak L}_2 \cdot \b, \
c_2 \in {\mathfrak L}_2 \cdot\b\setminus {\mathfrak L}_1\cdot \b$$ and put
$c=c_1+c_2$; clearly,
$$c \not\in {\mathfrak L}_1\cdot \b, \ c \not\in  {\mathfrak L}_2\cdot \b.$$
In all three cases
$$c\OC_F={\mathfrak M}\cdot \b$$
where the ideal ${\mathfrak M}=\prod_{\mathfrak P} {\mathfrak P}^{m_{\mathfrak
P}}$ is a (finite) product of powers of (non-zero) prime ideals ${\mathfrak
P}$, none of which divides $\ell$. It follows that $\b/c \OC_F$ is a (finite)
$\OC_F/{\mathfrak M}$-module. By the Chinese Remainder Theorem,
$$\OC_F/{\mathfrak M}=\oplus _{\mathfrak P}\OC_F/{\mathfrak P}^{m_{\mathfrak
P}}.$$ Therefore $\b/c \OC_F$ is a product of finite $\OC_F/{\mathfrak
P}^{m_{\mathfrak P}}$-modules. Since the multiplication by the residual
characteristic of ${\mathfrak P}$ kills $\OC_F/{\mathfrak P}$, it follows that
the ${m_{\mathfrak P}}$th power of this characteristic kills every
$\OC_F/{\mathfrak P}^{m_{\mathfrak P}}$-module. This implies that the order of
$\b/c \OC_F$ is a product of powers of residual characteristics of ${\mathfrak
P}$'s and therefore is not divisible by $\ell$.
\end{proof}

\begin{cor}
\label{degE} For every prime $\ell$ there exists an isogeny $E \to E^{(\b)}$,
whose degree is not divisible by $\ell$.
\end{cor}

\begin{sect}
{\bf The construction}. Choose an imaginary quadratic field $F$ with class
number $>1$ and pick a {\sl non}-principal ideal $\b\subset \OC_F$. We have
$$\j(E^{(\b)}) \ne \j(E).$$
 There exists an algebraic number field $L \subset \C$ such that:

\begin{itemize}
\item
  $L$
contains $F$, $\j(E)$ and $\j(E^{(\b)})$.

\item
 The elliptic curves $E$ and $E^{(\b)}$ are
defined over $L$.

\item All homomorphisms between $E$ and $E^{(\b)}$ are defined over $L$.
\end{itemize}

 Let us choose a maximal ideal $\q\subset \OC_F$ such that both $E$ and $E^{(\b)}$
have good reduction at $\q$ and  $j(E)- j(E^{(\b)})$ does {\sl not}
lie in $\q$. (Those conditions are satisfied by all but finitely
many $\q$.) Let $K$ be the (finite) residue field at $\q$, let ${\bf
E}$ and ${\bf E}^{(b)}$ be the reductions at $\q$ of $E$ and
$E^{(\b)}$ respectively: they are elliptic curves over $K$. Then
$\j({\bf E})$ and $\j({\bf E}^{(b)})$ are the reductions modulo $\q$
of $\j(E)$ and $\j(E^{(\b)})$ respectively. Our assumptions on $\q$
imply that
$$\j({\bf E}) \ne \j({\bf E}^{(b)}).$$
Therefore  ${\bf E}$ and ${\bf E}^{(b)}$ are not isomorphic over $K$
and even over $\bar{K}$!

On the other hand, it is known \cite[Ch. 9, Sect. 3]{LangE} that there is a
natural embedding
$$\Hom(E,E^{(\b)}) \hookrightarrow \Hom({\bf E},{\bf E}^{(\b)})$$
that respects the degrees of isogenies. It follows from Corollary \ref{degE}
that for every prime $\ell$ there exists an isogeny ${\bf E} \to {\bf
E}^{(\b)}$, whose degree is not divisible by $\ell$. Now Proposition
\ref{composite} implies that ${\bf E}_n \cong {{\bf E}^{(\b)}}_n$ for all
positive integers $n$. It follows from Corollary \ref{isomDIV} that the
$\ell$-divisible groups ${\bf E}(\ell)$ and ${\bf E}^{(\b)}(\ell)$ are
isomorphic for all $\ell$, including $\ell=\fchar(K)$. Since both ${\bf
E}(\bar{K})$ and ${\bf E}^{(\b)}(\bar{K})$ are torsion groups, they are
isomorphic as Galois modules. This implies that their subgroups of all Galois
invariants are isomorphic, i.e., the finite groups ${\bf E}({K})$ and ${\bf
E}^{(\b)}({K})$ are isomorphic.
\end{sect}

\end{document}